\documentclass[oneside,english]{amsart}
\usepackage[T1]{fontenc}
\usepackage[latin9]{inputenc}
\usepackage{geometry}
\usepackage{amstext}
\usepackage{amsthm}
\usepackage{amssymb}

\makeatletter
\numberwithin{equation}{section}
\numberwithin{figure}{section}
  \theoremstyle{remark}
  \newtheorem*{rem*}{\protect\remarkname}
  \theoremstyle{remark}
  \newtheorem*{acknowledgement*}{\protect\acknowledgementname}

\AtBeginDocument{
  
}

\makeatother

\usepackage{babel}
  \providecommand{\acknowledgementname}{Acknowledgement}
  \providecommand{\remarkname}{Remark}

\begin{document}

\title{An example of non-embeddability of the Ricci flow}

\author{Mohammad Safdari}
\begin{abstract}
For an evolution of metrics $(M,g_{t})$ there is a $t$-smooth family
of embeddings $e_{t}:M\rightarrow\mathbb{R}^{N}$ inducing $g_{t}$,
but in general there is no family of embeddings extending a given
initial embedding $e_{0}$. We give an example of this phenomenon
when $g_{t}$ is the evolution of $g_{0}$ under the Ricci flow. We
show that there are embeddings $e_{0}$ inducing $g_{0}$ which do
not admit of $t$-smooth extensions to $e_{t}$ inducing $g_{t}$
for any $t>0$. We also find hypersurfaces of $\textrm{dim}>2$ that
will not remain a hypersurface under Ricci flow for any positive time.\thanks{$^{1}\;$Department of Mathematical Sciences, Sharif University of
Technology, Tehran, Iran\protect \\
Email address: safdari@sharif.edu} 
\end{abstract}

\maketitle

\section{Preliminaries}

Let $f:\mathbb{R}^{n}\rightarrow\mathbb{R}$ be a smooth function
and let
\[
M^{n}=\textrm{graph}(f)\hookrightarrow\mathbb{R}^{n+1}
\]
Let $\langle\;,\;\rangle$ be the standard inner product on $\mathbb{R}^{n+1}$,
and $g$ be the induced metric on $M$ from $\mathbb{R}^{n+1}$. Let
$\nabla\,,\,D$ be the standard covariant derivatives on $\mathbb{R}^{n}\,,\,\mathbb{R}^{n+1}$
respectively. First we will compute the metric and curvature of $M$.
Note that $M$ is diffeomorphic to $\mathbb{R}^{n}$ and we can cover
it by one chart, which we will do from now on. Now 
\[
F=(\textrm{id},f)\;:\;M\rightarrow\mathbb{R}^{n+1}
\]
is an embedding, so the tangent vectors to $M$ are $\partial_{i}F=(e_{i},\partial_{i}f)$
where $e_{i}$s are the standard basis of $\mathbb{R}^{n}$. The
components of $g$ are 
\begin{equation}
g_{ij}=\langle\partial_{i}F,\partial_{j}F\rangle=\delta_{ij}+\partial_{i}f\partial_{j}f
\end{equation}
The unit normal to $M$ is 
\begin{equation}
N=\frac{1}{\sqrt{1+|\overset{}{\nabla}f|^{2}}}(\overset{}{\nabla}f,-1)
\end{equation}
Also the components of the second fundamental form are 
\[
h_{ij}=-\langle D_{\partial_{i}F}\partial_{j}F,N\rangle
\]
Since
\[
D_{\partial_{i}F}\partial_{j}F=\underset{k}{\sum}\partial_{i}F^{k}\partial_{k}\partial_{j}F=\partial_{i}\partial_{j}F+\partial_{i}f\partial_{n+1}\partial_{j}F=\partial_{i}\partial_{j}F=(0,\partial_{i}\partial_{j}f)
\]
(because $\partial_{j}F$ is independent of $x_{n+1}$) we have
\begin{equation}
h_{ij}=-\langle\partial_{i}\partial_{j}F,N\rangle=\frac{\partial_{i}\partial_{j}f}{\sqrt{1+|\overset{}{\nabla}f|^{2}}}
\end{equation}
The Gauss equation implies
\begin{equation}
R_{ijkl}=h_{il}h_{jk}-h_{ik}h_{jl}=\frac{1}{1+|\nabla f|^{2}}(\partial_{i}\partial_{l}f\,\partial_{j}\partial_{k}f-\partial_{i}\partial_{k}f\,\partial_{j}\partial_{l}f)
\end{equation}
By an easy induction on $n$ we find
\begin{equation}
\textrm{det}\,g_{ij}=1+\sum(\partial_{i}f)^{2}=1+|\nabla f|^{2}
\end{equation}
and also the components of $g^{-1}$
\begin{equation}
g^{ij}=\delta_{ij}-\frac{\partial_{i}f\partial_{j}f}{1+|\nabla f|^{2}}
\end{equation}
Now we can compute the Christoffel symbols
\begin{equation}
\Gamma_{ij}^{k}=\frac{1}{2}g^{kl}(\partial_{i}g_{jl}+\partial_{j}g_{il}-\partial_{l}g_{ij})=\frac{\partial_{k}f\partial_{i}\partial_{j}f}{1+|\nabla f|^{2}}
\end{equation}

\section{The Example}

Let $f:\mathbb{R}^{n}\rightarrow\mathbb{R}$ be given by 
\begin{equation}
f(x_{1},\dots,x_{n})=\underset{r,q}{\sum}a_{rq}x_{r}x_{q}^{2}
\end{equation}
(where $(a_{rq})$ is not necessarily symmetric) We are interested
in the evolution of $(M,g)$ under the Ricci flow. Let $p$ denote
the origin in $\mathbb{R}^{n+1}$. The derivatives of $f$ are
\begin{equation}
\partial_{i}f=\underset{r\neq i}{\sum}a_{ri}x_{r}x_{i}+\underset{q\neq i}{\sum}a_{iq}x_{q}^{2}+3a_{ii}x_{i}^{2}
\end{equation}

\begin{equation}
\partial_{i}\partial_{j}f=\begin{cases}
2a_{ij}x_{j}+2a_{ji}x_{i} & i\ne j\\
\\
\underset{q\neq i}{\sum}2a_{qi}x_{q}+6a_{ii}x_{i} & i=j
\end{cases}
\end{equation}
Note that all these expressions vanish at the origin, so both the
curvature and the connection vanish at $p$. In addition we have $g_{ij}=\delta_{ij}$
at $p$. 

We know that under the Ricci flow the Riemann curvature tensor evolves
as 
\[
\frac{\partial}{\partial t}Rm=\triangle Rm+Rm*Rm+Rm*Ric
\]
where $A*B$ is a sum of contractions of components of the tensors
$A$ and $B$ by the metric. Now if we look at this equation at $x=p$
and $t=0$ we get $ $ 
\begin{equation}
\frac{\partial}{\partial t}\bigg|_{t=0}Rm=\triangle Rm
\end{equation}
But 
\[
\triangle Rm=g^{ij}\partial_{i}\partial_{j}Rm+\Gamma*\partial Rm+\partial\Gamma*Rm+\Gamma*Rm
\]
Since $\Gamma$ and $Rm$ vanish at $p$ and $g$ is the identity
there, we obtain$\frac{}{}$ 
\begin{equation}
\frac{\partial}{\partial t}\bigg|_{t=0}Rm=\sum\partial_{i}\partial_{i}Rm
\end{equation}
at $x=p$. We define the tensor $A$ by 
\begin{equation}
A_{ijkl}:=\partial_{i}\partial_{l}f\,\partial_{j}\partial_{k}f-\partial_{i}\partial_{k}f\,\partial_{j}\partial_{l}f
\end{equation}
then $Rm=\frac{1}{1+|\nabla f|^{2}}\,A$, and as $A,\,\partial A$
vanish at $p$, we have 
\begin{equation}
\partial\partial Rm=\frac{1}{1+|\nabla f|^{2}}\,\partial\partial A
\end{equation}
at $p$. Since $A$ has the symmetries of the curvature tensor, it
is enough to compute $A_{ijkl}$ for $i<j,k<l,i\leq k$
\[
A_{ijkl}=\begin{cases}
4(a_{il}x_{l}+a_{li}x_{i})(a_{kj}x_{j}+a_{jk}x_{k})\\
-4(a_{ik}x_{k}+a_{ki}x_{i})(a_{jl}x_{l}+a_{lj}x_{j}) & i,j,k,l\;\textrm{are all distinct}\\
\\
4(a_{il}x_{l}+a_{li}x_{i})(\underset{q\neq j}{\sum}a_{qj}x_{q}+3a_{jj}x_{j})\\
-4(a_{ij}x_{j}+a_{ji}x_{i})(a_{jl}x_{l}+a_{lj}x_{j}) & i<k=j<l\\
\\
4(a_{ij}x_{j}+a_{ji}x_{i})(a_{kj}x_{j}+a_{jk}x_{k})\\
-4(a_{ik}x_{k}+a_{ki}x_{i})(\underset{q\neq j}{\sum}a_{qj}x_{q}+3a_{jj}x_{j}) & i<k<l=j\\
\\
4(a_{il}x_{l}+a_{li}x_{i})(a_{ij}x_{j}+a_{ji}x_{i})\\
-4(a_{jl}x_{l}+a_{lj}x_{j})(\underset{q\neq i}{\sum}a_{qi}x_{q}+3a_{ii}x_{i}) & i=k,j\neq l\\
\\
4(a_{ij}x_{j}+a_{ji}x_{i})^{2}\\
-4(\underset{q\neq i}{\sum}a_{qi}x_{q}+3a_{ii}x_{i})(\underset{r\neq j}{\sum}a_{rj}x_{r}+3a_{jj}x_{j}) & i=k,j=l
\end{cases}
\]
Therefore 
\[
\sum\partial_{s}\partial_{s}A_{ijkl}=\begin{cases}
0 & i,j,k,l\;\;\\
 & \textrm{are all distinct}\\
8(a_{il}a_{lj}+a_{li}a_{ij})-8a_{ij}a_{lj} & i<k=j<l\\
\\
8a_{ij}a_{kj}-8(a_{ik}a_{kj}+a_{ki}a_{ij}) & i<k<l=j\\
\\
8a_{li}a_{ji}-8(a_{jl}a_{li}+a_{lj}a_{ji}) & i=k,j\neq l\\
\\
8(a_{ij}^{2}+a_{ji}^{2})-8(\underset{q\neq i,j}{\sum}a_{qi}a_{qj}+3a_{ii}a_{ij}+3a_{jj}a_{ji}) & i=k,j=l
\end{cases}
\]
Now if we choose $(a_{rq})$ such that 
\begin{equation}
a_{\alpha\beta}a_{\beta\gamma}+a_{\beta\alpha}a_{\alpha\gamma}=a_{\alpha\gamma}a_{\beta\gamma}\label{* eq}
\end{equation}
for $\alpha,\beta,\gamma$ mutually distinct then all non-diagonal
entries of $\sum\partial_{s}\partial_{s}A_{ijkl}$ and hence $\frac{\partial}{\partial t}Rm$
vanish at $x=p$ and $t=0$ (note that $\frac{1}{1+|\nabla f|^{2}}=1$
at $p$). 

Let $a_{ij}=0$ for $i<j$ and $a_{ji}=1$, then (\ref{* eq}) holds.
Also let the diagonal elements $a_{ii}$ be $1$, then the diagonal
entries of $\frac{\partial}{\partial t}Rm(0,p)$ are all negative.
Therefore as $Rm(0,p)=0$, the sectional curvatures of $Rm(t,p)$
will be negative for small $t$. Therefore for any positive $t$,
$(M,g(t))$ is no longer a hypersurface when $n\geq3$. Note that
no neighborhood of $p$ can be embedded in $\mathbb{R}^{n+1}$ for
any $t>0$ and using this we can construct closed hypersurfaces that
will not remain a hypersurface for any positive time under Ricci flow.

We also observe that for $n\geq2$ if we consider the embedding 
\[
\varphi:M^{n}\hookrightarrow\mathbb{R}^{n+1}\hookrightarrow\mathbb{R}^{n+k}
\]
for any $k\geq1$, then there is no evolution of $\varphi$ that induces
the Ricci flow on $M$. In fact, if such evolution of $\varphi$ exists,
then by the Gauss equation we will have
\[
Rm(t)(X,Y,Z,W)=\langle\Pi(t)(X,W),\Pi(t)(Y,Z)\rangle-\langle\Pi(t)(X,Z),\Pi(t)(Y,W)\rangle
\]
where $\Pi$ is the second fundamental form. Differentiating we obtain
\begin{eqnarray*}
\frac{\partial}{\partial t}Rm(0)(X,Y,Z,W) &  & =\langle\frac{\partial}{\partial t}\Pi(0)(X,W),\Pi(0)(Y,Z)\rangle\\
 &  & +\,\langle\Pi(0)(X,W),\frac{\partial}{\partial t}\Pi(0)(Y,Z)\rangle\\
 &  & -\,\langle\frac{\partial}{\partial t}\Pi(0)(X,Z),\Pi(0)(Y,W)\rangle\\
 &  & -\,\langle\Pi(0)(X,Z),\frac{\partial}{\partial t}\Pi(0)(Y,W)\rangle
\end{eqnarray*}
at $p$. But $\Pi(0)=0$ at $p$ and this contradicts the fact that
$\frac{\partial}{\partial t}Rm(0)\neq0$.

Also note that even changing the metric on $\mathbb{R}^{n+k}$ will
not allow the existence of an evolution of $\varphi$ that induces
$g(t)$ since in this case
\[
Rm(t)(X,Y,Z,W)=\eta(t)(\Pi(t)(X,W),\Pi(t)(Y,Z))-\eta(t)(\Pi(t)(X,Z),\Pi(t)(Y,W))
\]
where $\eta(t)$ is the evolution of the standard metric on $\mathbb{R}^{n+k}$.
Thus
\begin{eqnarray*}
\frac{\partial}{\partial t}Rm(0)(X,Y,Z,W) &  & =\langle\frac{\partial}{\partial t}\Pi(0)(X,W),\Pi(0)(Y,Z)\rangle\\
 &  & +\,\langle\Pi(0)(X,W),\frac{\partial}{\partial t}\Pi(0)(Y,Z)\rangle\\
 &  & -\,\langle\frac{\partial}{\partial t}\Pi(0)(X,Z),\Pi(0)(Y,W)\rangle\\
 &  & -\,\langle\Pi(0)(X,Z),\frac{\partial}{\partial t}\Pi(0)(Y,W)\rangle\\
 &  & +\,\frac{\partial}{\partial t}\eta(0)(\Pi(0)(X,Z),\Pi(0)(Y,W))
\end{eqnarray*}
at $p$. Again $\Pi(0)=0$ at $p$ and we get a contradiction with
$\frac{\partial}{\partial t}Rm(0)\neq0$.
\begin{rem*}
For a generic isometric embedding of a Riemannian manifold $(M,g)$
in $\mathbb{R}^{N}$, the metric $g(t)$ can be embedded in $\mathbb{R}^{N}$
for small $t>0$. The problem is that we do not know if and when the
evolution of the metric encounters an obstacle beyond which we may
not be able to extend the embedding. A successful resolution of this
problem will have interesting consequences. For example it may allow
us to obtain isometric embeddings of a surface of $\textrm{genus}>1$
with constant negative curvature in $\mathbb{R}^{5}$.
\end{rem*}
\begin{acknowledgement*}
The author would like to thank Mehrdad Shahshahani and Burkhard Wilking
for their help with this research.
\end{acknowledgement*}
\bibliographystyle{plain}
\nocite{*}
\bibliography{ref}

\end{document}